\documentclass{article}

\usepackage[backend=bibtex, 
sortcites=true,
firstinits=true,
hyperref,
maxbibnames=99,
]{biblatex}

\bibliography{Bibliography}{}

\usepackage{enumerate}

\usepackage[centertags]{amsmath}
\usepackage{amsfonts}
\usepackage{amssymb}
\usepackage{amsthm}
\usepackage{newlfont}
\usepackage{mathtools}
\usepackage[top=1in, bottom=1in, left=1in, right=1in]{geometry}
\usepackage{tikz}

\theoremstyle{plain}
\newtheorem{theorem}{Theorem}[subsection]
\newtheorem{proposition}[theorem]{Proposition}
\newtheorem{corollary}[theorem]{Corollary}
\newtheorem{lemma}[theorem]{Lemma}

\theoremstyle{definition}
\newtheorem{definition}[theorem]{Definition}

\newcommand{\R}{\mathbb{R}}
\newcommand{\N}{\mathbb{N}}
\newcommand{\Z}{\mathbb{Z}}

\newcommand{\Deg}{\operatorname{Deg}}
\newcommand{\cost}{\operatorname{cost}}
\newcommand{\diag}{\operatorname{diag}}
\newcommand{\diam}{\operatorname{diam}}
\newcommand{\sgn}{\operatorname{sgn}}

\newcommand{\Hidden}[1]{}

\begin{document}

\title{Non-negative Ollivier curvature on graphs, reverse Poincar\'e inequality, Buser inequality, Liouville property, Harnack inequality and eigenvalue estimates}
\author{
%Marzieh Eidi\footnote{MPI Leipzig, eidi@mis.mpg.de}~~~~~
%~~~~~J\"urgen Jost\footnote{MPI Leipzig, jost@mis.mpg.de}~~~~~
%~~~~~
Florentin M\"unch\footnote{MPI MiS Leipzig, muench@mis.mpg.de}%~~~~~Christian Rose\footnote{MPI Leipzig, crose@mis.mpg.de}
}
\date{\today}
\maketitle

\begin{abstract}
%We prove the reverse Poincar\'e inequality for graphs with non-negative Ollivier curvature. 
We prove that for combinatorial graphs with non-negative Ollivier curvature, one has
\[
\|P_t \mu - P_t \nu\|_1   \leq   \frac{W_1(\mu,\nu)}{\sqrt{t}}
\]
for all probability measures $\mu,\nu$ where $P_t$ is the heat semigroup and $W_1$ is the $\ell_1$-Wasserstein distance.
This turns out to be an equivalent formulation of a version of reverse Poincar\'e inequality.
Furthermore, this estimate allows us to prove Buser inequality, Liouville property and the the eigenvalue estimate $\lambda_1 \geq \log(2)/\diam^2$.
\end{abstract}

%%%%%%%%%%%%%%
%\section{Introduction}
%%%%%%%%%%%%%%
\tableofcontents

\section{Introduction}
Discrete Ricci curvature is a vibrant subject of research in the recent years. Particularly interesting is the case of non-negative Ricci curvature where various fundamental results are known for Riemannian manifolds such as Li-Yau inequality, Harnack inequality, Buser inequality and eigenvalue estimates.
Great effort has been invested to prove discrete analogs of these results \cite{bauer2015li,munch2014li,
liu2014eigenvalue,liu2015curvature,
liu2018buser,erbar2018poincare,
chung2014harnack,horn2014volume}.
While Ricci curvature in Riemannian geometry is a well known notion since over a century,  various non-equivalent notions of discrete Ricci curvature have been developed in the last decade.
Particularly, we mention Bakry-Emery curvature based on Bochner's formula \cite{schmuckenschlager1998curvature,lin2010ricci}, entropic curvature based on the convexity of the entropy of transport geodesics \cite{erbar2012ricci}, and Ollivier curvature based on optimal transport \cite{ollivier2007ricci,ollivier2009ricci}.
Although there does not seem to hold any implication between these three Ricci curvature notions, every of these Ricci curvature notions allows for a characterization of  non-negative curvature via a gradient estimate of the form
\[
\|\nabla P_t f\| \leq \|\nabla f\|_{P_t}
\]
when choosing appropriate notions of the norms $\|\cdot \|$ and $\|\cdot\|_{P_t}$ where $P_t$ is the heat semigroup.
It turned out that for proving Harnack inequality, Buser inequality and eigenvalue estimates, the key ingredient is the reverse Poincar\'e inequality stating that for all  bounded functions $f$, one has
\[
\|\nabla P_t f\|_\infty \leq const. \cdot \frac{\|f\|_\infty}{\sqrt {t}}
\]
where $\|\nabla\cdot\|_\infty$ is the optimal Lipschitz constant w.r.t the combinatorial distance.
While Ledoux's method provides an elegant way to prove the reverse Poncar\'e inequality from the gradient estimate in case of non-negative Bakry Emery curvature and in case of non-negative entropic Ricci curvature, this method seems to fail in the case of Ollivier curvature. Even more remarkable, the gradient estimate
\[
\|\nabla P_t f\|_\infty \leq \|f\|_\infty
\]
and its dual version, the Wasserstein contraction property, seem to be the only investigated analytic implications of non-negative Ollivier curvature, up to the recently established Liouville property \cite{jost2019Liouville} which will also be improved in this paper.

In contrast to non-negative curvature, the case of a positive lower curvature bounds is well studied for both, Bakry Emery and Ollivier curvature, and results in analogy to Lichnerowicz eigenvalue estimate, Bonnet-Myers diameter bounds with Cheng's rigidity theorem, Gaussian concentration and volume growth bounds have been proven \cite{cushing2018rigidity,liu2017rigidity,
ollivier2007ricci,ollivier2009ricci,
bhattacharya2015exact,chung2017curvature,
fathi2015curvature,jost2019Liouville,
jost2014ollivier,liu2016bakry,paeng2012volume}.
Example classes of graphs with positive curvature have been studied in \cite{cushing2016bakry,yamada2017curvature,
bonini2019condensed,cushing2018curvature}

The main goal of this paper is to prove the reverse Poincar\'e inequality under non-negative Ollivier curvature (see Theorem~\ref{thm:NonNegOllRevPoin}) by applying the coupling method and by showing the existence of suitable transport plans (see Theorem~\ref{thm:transportMass}) by which we can compare the coupled heat equation to the heat equation on an unweighted birth death chain (see Theorem~\ref{thm:semigroupCouplingEstimate}).
This result is of particular relevance since Ollivier curvature seems to be the most prominent notion of discrete Ricci curvature in terms of network analysis \cite{sandhu2015analytical,
sandhu2015graph,sandhu2016ricci,
tannenbaum2015ricci,wang2014wireless,
wang2016interference,ni2015ricci,
ni2018network,ni2019community}.

The reverse Poincar\'e inequality turns out to be equivalent to the total variation estimate
\[
\|P_t \mu - P_t \nu\|_1   \leq   \frac{W_1(\mu,\nu)}{\sqrt{t}}
\]
for all probability measures $\mu,\nu$ where 
%$P_t$ is the heat semigroup and 
$W_1$ is the $\ell_1$-Wasserstein distance (see Theorem~\ref{thm:RevPoincarChar}).
As applications, we prove Buser inequality (see Theorem~\ref{thm:Buser}), Liouville property, Harnack inequality and eigenvalue estimates.

The paper is structured as follows:
In section \ref{s:Setup}, we introduce our general setup and notation including Ollivier and Bakry Emery curvature.
In section~\ref{s:CouplingMethod}, we apply a purely analytic version of the coupling method to prove the reverse Poincar\'e inequality under non-negative Ollivier curvature.
In section~\ref{s:Applications}, we give applications of the reverse Poincar\'e inequality. In particular, we characterize the reverse Poincar\'e via a total variation estimate, we prove Buser inequality, Liouville property, Harnack inequality and eigenvalue estimates.
Finally in section ~\ref{s:Discussion}, we give an informal random walk interpretation of the coupling method, and we discuss the fact that non-negative Ollivier curvature and non-negative Bakry Emery curvature imply precisely the same reverse Poincar\'e inequality.

\subsection{Setup and notation}\label{s:Setup}
A \emph{graph} $G=(V,q)$ consists of a countable set $V$ and a function $q:V\times V \to [0,\infty)$ such that $\#\{y: q(x,y)>0\} < \infty$
for all $x \in V$.
We say, $G$ is \emph{reversible} or \emph{undirected} if there exists $m:V \to (0,\infty)$ s.t.
\[
q(x,y)m(x) = q(y,x)m(y)
\]
for all $x,y \in V$.
The notion of reversibility comes from Markov chain theory and the notion of undirectedness comes from graph theory which is closely related to Markov chain theory. 
We will use several function spaces. These are $C(V):=\{f:V\to \R\}$ and $\ell_\infty(V):=\{f \in C(V): \|f\|_\infty = \sup_x |f(x)| < \infty\}$ and $\ell_1(V,m):=\{f \in C(V): \|f\|_1=\sum_x m(x)|f(x)| < \infty\}$.
The Laplacian $\Delta:C(V) \to C(V)$ is given by
\[
\Delta f(x) := \sum_y q(x,y)(f(y)-f(x)).
\]
In the notion of continuous time Markov chains, the function $q$ can be interpreted as transition  rate.
We write
\[
q_{\min}:= \inf \{q(x,y):x,y\in V, q(x,y)>0\}.
\]
Reversibility of the graph means that the Laplacian is a symmetric operator w.r.t the invariant measure $m$.
For reversible graphs, we write $x\sim y$ if $q(x,y)>0$ and we define the combinatorial distance
\[
d(x,y):=\inf\{n:x=x_0\sim \ldots \sim x_n=y\}.
\]
We now give the definition for Ollivier curvature which has been introduced in \cite{ollivier2009ricci, ollivier2007ricci} and generalized in \cite{jost2014ollivier,lin2011ricci, munch2017ollivier,eidi2019ollivier,
asoodeh2018curvature}.
For a reversible graph $G=(V,q)$ and $f \in C(V)$, and $x\neq y \in V$, we define
\[
\nabla_{xy} f := \frac{|f(x)-f(y)|}{d(x,y)}
\]
and
\[
\|\nabla f\|_\infty := \sup_{x\neq y} \nabla_{xy}f
\]
and
\[
\kappa(x,y) := \inf_{\substack{\nabla_{yx}f=1\\ \|\nabla f\|_\infty = 1}} \nabla_{xy}\Delta f.
\]
By \cite[Theorem~2.1]{munch2017ollivier}, the term $\kappa(x,y)$ is a generalization of Lin-Lu-Yau's modification of Ollivier curvature.
We remark that this modification corresponds to lazy random walks and is always larger or equal to the curvature corresponding to non-lazy random walks. Thus, non-negative curvature in our sense is a weaker assumption than non-negative curvature for non-lazy random walks.

We now give the definition of Bakry Emery curvature which goes back to \cite{BakryEmery85} for Riemannian geometry and has been introduced on discrete spaces in \cite{schmuckenschlager1998curvature} and later in \cite{lin2010ricci} independently. For $f,g \in C(V)$, let $\Gamma_k:C(V)\times C(V) \to C(V)$ be inductively given by $\Gamma_0(f,g):=f\cdot g$ and
\[
2\Gamma_{k+1}(f,g) := \Delta \Gamma_k(f,g) - \Gamma_k(f,\Delta g) - \Gamma_k(\Delta f,g).
\]
We write $\Gamma_k f := \Gamma_k(f,f)$ and $\Gamma:=\Gamma_1$. We say a graph satisfies the \emph{curvature dimension condition} $CD(K,n)$ for $K\in \R$ and $n \in (0,\infty]$ if
\[
\Gamma_2 f \geq K\Gamma f +\frac 1 n (\Delta f)^2
\]
for all $f \in C(V)$. We say a graph has non-negative Bakry Emery curvature if it satisfies $CD(0,\infty)$.

\section{Coupling method and non-negative Ollivier curvature}\label{s:CouplingMethod}
The coupling method is a powerful theory which has been invented to describe couplings of two random walks. In this paper, we give a purely analytic description of the coupling method by  using the heat semigroup instead of random walks which are in one-to-one correspondence due to Feynman-Kac formula.

\subsection{Optimal transport}

Due to Kantorovich duality, the Ollivier curvature can also be calculated via transport plans, see Proposition~\ref{pro:CharTransport}.

Let $G=(V,q)$ be a graph and let $x_0,y_0 \in V$.
A map $\rho = \rho_{x_0y_0} : V\times V \to [0,\infty)$ is called called a \emph{transport plan} from $x_0$ to $y_0$ if
\begin{align}
\sum_{y \in V} \rho(x,y) &= q(x_0,x)  \qquad \mbox{ for all } x \in V \setminus \{x_0\} \mbox{ and} \label{eq:rhoXProp}\\
\sum_{x \in V} \rho(x,y) &= q(y_0,y) \qquad \mbox{ for all } y \in V \setminus \{y_0\}  \label{eq:rhoYProp}.
\end{align}
We remark that due to local finiteness of $G$, every transport plan is finitely supported.
The \emph{transport cost} of a transport plan $\rho$ is given by
\[
\cost(\rho) := \sum_{x,y \in V}\rho(x,y) \big(d(x_0,y_0) - {d(x,y)}\big).
\]
The mass $\mu_\rho(k)$ which is transported by $\rho$ over a distance $d(x_0,y_0)+k$ can be calculated via
\[
\mu_\rho(k) := \sum_{\substack{x,y \in V \\ d(x,y)- d(x_0,y_0)=k}}\rho(x,y).
\]
Thus,
\[
\cost(\rho) = - \sum_k k\mu_\rho(k).
\]
A transport plan which maximizes the cost is called optimal.
Properties of optimal transport plans have been widely used to investigate Ollivier curvature
\cite{bhattacharya2015exact,bourne2017ollivier,
loisel2014ricci,lin2011ricci,jost2014ollivier}.
Specifically due to Kantorovic duality, the Ollivier curvature can also be expressed via the cost of a transport plan.
\begin{proposition}[Curvature via transport cost, see {{\cite[Proposition~2.4]{munch2017ollivier}}}]\label{pro:CharTransport}
Let $G=(V,w,m)$ be a graph and let $x_0 \neq y_0$ be vertices. Then,
\begin{align}
d(x_0,y_0)\kappa(x_0,y_0) &= \sup_{\rho}  \cost(\rho) \label{eq:PropTransport} 
\end{align}
where the supremum is taken over all transport plans $\rho$.
\end{proposition}

We now give a new result concerning optimal transport plans stating that there always exists an optimal transport plan between $x_0\neq y_0 \in V$ such that a significant amount of mass is transported over a distance $d(x_0,y_0)-1$, and such that no mass is transported over the distance $d(x_0,y_0) \pm 2$.
\begin{theorem}\label{thm:transportMass}
Let $G=(V,q)$ be a reversible graph.

Let $x_0 \neq y_0 \in V$.
Then, there exists an optimal transport plan $\rho$ from $x_0$ to $y_0$ s.t.
\begin{itemize}
\item $\mu_\rho(k)=0$ whenever $|k|>1$,
\item $\mu_\rho(-1) \geq 2 q_{\min}$. 
\end{itemize}
Moreover every optimal transport plan $\rho$ from $x_0$ to itself satisfies
\[
\mu_\rho(k)=0
\]
for all $k \neq 0$.
\end{theorem}
\begin{proof}
We first prove the 'moreover' statement.
A transport plan form $x_0$ to itself is given by $\rho_0(x,x) = q(x_0,x)$ and $\rho_0(x,y)=0$ if $x\neq 0$. Obviously, $\cost(\rho_0) = 0$.
For every transport plan $\rho$ from $x_0$ to itself, we have $\mu_\rho(k)=0$ if $k < 0$. Thus if $\rho$ is optimal, then $\cost(\rho) = 0$ and hence, $\mu_\rho(k)=0$ for $k > 0$ which finishes the prove of the 'moreover' statement.

We now prove that for $x_0 \neq y_0$, there exists a transport plan $\rho$ with the desired properties.
%Let $x_0\neq y_0 \in V$.
Let $\rho_0$ be an optimal transport plan from $x_0$ to $y_0$.
We aim to construct $\pi:V\times V\to \R$ s.t.
$\rho = \rho_0 + \pi$ is an optimal transport plan with the desired properties.
We aim to find an appropriate set $A \subset V^2$ where we force $\rho$ to vanish. The set $A$ will be the union of the following sets.
We define
\begin{align*}
H &:= \{(x,y)\in B_1(x_0) \times B_1(y_0) : d(x,y)-d(x_0,y_0) \geq 2 \},\\
L &:= \{(x,y)\in B_1(x_0) \times B_1(y_0) : d(x,y)-d(x_0,y_0) \leq -2 \},\\
X&:= \{(x,y)\in B_1(x_0) \times B_1(y_0) :d(x,y_0) < d(x_0,y_0) \mbox{ and } d(x,y) \geq d(x_0,y_0) \},\\
Y&:= \{(x,y)\in B_1(x_0) \times B_1(y_0) : d(x_0,y) < d(x_0,y_0) \mbox{ and } d(x,y) \geq d(x_0,y_0)\},\\
A &:= L\cup H \cup X \cup Y.
\end{align*}
Observe that this union is disjoint.
We will show that $\rho$ vanishing on $H\cup L$ gives $\mu_\rho(k)=0$ whenever $|k|>1$, and that $\rho$ vanishing on $X\cup Y$ gives $\mu_{\rho}(-1) + \mu_{\rho}(-2)\geq 2q_{\min}$.
We start collecting basic properties of the sets defined above.
By triangular inequality for all $(x,y)\in V\times V$, 
\begin{align*}
(x,y) \in  H \qquad &\Longrightarrow \qquad d(x,y)-d(x_0,y_0) = 2,\\
(x,y) \in  L \qquad &\Longrightarrow \qquad d(x,y)- d(x_0,y_0)=-2,\\
(x,y)  \in  X \cup Y \qquad &\Longrightarrow \qquad d(x,y)-d(x_0,y_0)=0.
\end{align*}
In particular,
\begin{align}\label{eq:proofTransEvenDist}
d(x,y) \in \{-2,0,2\} \qquad \mbox{ for all } (x,y) \in A.
\end{align}
Furthermore, $A \subseteq S_1(x_0) \times S_1(y_0)$.
For $x,y \in V$, we define
\[
\pi(x,y) := \begin{cases}
-\rho_0(x,y) &: (x,y) \in A, \\
\sum_{(x,\widetilde y) \in A} \rho_0(x,\widetilde y) &: y=y_0, \\
\sum_{(\widetilde x, y) \in A} \rho_0(\widetilde x, y) &: x=x_0, \\
0&:\mbox{else.}
\end{cases}
\]
Step 1: We show that $\rho=\rho_0 + \pi$ is a transport plan. We remark that the definition of $\cost(\pi)$ and $\mu_\pi(k)$ transfer canonically although $\pi$  not  a transport plan.
For $x \neq x_0$,
\begin{align*}
\sum_y \pi(x,y) &= \sum_{y\neq y_0} \pi(x,y) + \pi(x,y_0)\\
& = -\sum_{\substack{y\neq y_0\\(x,y)\in A}} \rho_0(x,y) +  \sum_{(x,\widetilde y) \in A} \rho_0(x,\widetilde y) \\
&=0
\end{align*}
where the last identity follows since $A \subseteq S_1(x_0) \times S_1(y_0)$. Similarly, $\sum_x \pi(x,y)=0$ for $y \neq y_0$.
Moreover, $\rho:=\rho_0 + \pi \geq 0$ and thus, $\rho$ is also a transport plan.
%We observe $\rho(x,y)=0$ for all $(x,y) \in A$.

Step 2: We show that $\rho = \rho_0 + \pi$ is optimal.
To this end, we have to show $\cost(\pi) \geq 0$.
For $W \subset V^2$, we write 
\[
\rho_0(W) := \sum_{(x,y)\in W} \rho_0(x,y).
\]
We have
\[
\mu_\pi(2) = -\sum_{\substack{(x,y)\in A\\d(x,y) - d(x_0,y_0) = 2} } \rho_0(x,y) = -\rho_0(H)
\]
and similarly,
\[
\mu_\pi(-2) = -\sum_{\substack{(x,y)\in A\\d(x,y) - d(x_0,y_0) = -2} } \rho_0(x,y) = -\rho_0(L)
\]

Subsequently, we will use that for all $(x,y)\in V^2$,
\begin{align}\label{eq:proofTransportHLimpl}
\begin{aligned}
(x,y) \in H  \quad \Longrightarrow  \quad d(x,y_0)=d(x_0,y) > d(x_0,y_0),\\
(x,y) \in L  \quad \Longrightarrow \quad d(x,y_0)=d(x_0,y) < d(x_0,y_0).
\end{aligned}
\end{align}
Thus,
\begin{align*}
\mu_\pi(1) 
&
=  \sum_{d(x,y_0) > d(x_0,y_0)} \pi(x,y_0) + 
\sum_{d(x_0,y) > d(x_0,y_0)} \pi(x_0,y) -\sum_{\substack{(x,y)\in A \\ d(x,y)-d(x_0,y_0)=1}} \rho_0(x,y) \\
&
\stackrel{\eqref{eq:proofTransEvenDist}}{=} \sum_{d(x,y_0) > d(x_0,y_0)} \pi(x,y_0) + 
\sum_{d(x_0,y) > d(x_0,y_0)} \pi(x_0,y) \\
&
= \sum_{\substack{(x,y)\in A \\ d(x_0,y) > d(x_0,y_0)}}\rho_0(x,y) + \sum_{\substack{(x,y)\in A \\ d(x,y_0) > d(x_0,y_0)}}\rho_0(x,y)\\
&
\stackrel{\eqref{eq:proofTransportHLimpl}}{=}
2{\rho_0}(H) + \sum_{\substack{(x,y)\in X\cup Y \\ d(x_0,y) > d(x_0,y_0)}} \rho_0(x,y) + \sum_{\substack{(x,y)\in X\cup Y \\ d(x,y_0) > d(x_0,y_0)}} \rho_0(x,y) \\
&
=
2{\rho_0}(H) + \sum_{\substack{(x,y)\in X \\ d(x_0,y) > d(x_0,y_0)}} \rho_0(x,y) + \sum_{\substack{(x,y)\in Y \\ d(x,y_0) > d(x_0,y_0)}} \rho_0(x,y) \\
&
\leq 2\rho_0(H) + \rho_0(X) + \rho_0(Y).
\end{align*}
Similarly,
\begin{align*}
\mu_\pi(-1) &=  \sum_{d(x,y_0) < d(x_0,y_0)} \pi(x,y_0) + 
\sum_{d(x_0,y) < d(x_0,y_0)} \pi(x_0,y) -\sum_{\substack{(x,y)\in A \\ d(x,y)-d(x_0,y_0)=-1}} \rho_0(x,y)  \\
&\stackrel{\eqref{eq:proofTransEvenDist}}{=} \sum_{d(x,y_0) < d(x_0,y_0)} \pi(x,y_0) + 
\sum_{d(x_0,y) < d(x_0,y_0)} \pi(x_0,y)\\
&= \sum_{\substack{(x,y)\in A \\ d(x_0,y) < d(x_0,y_0)}}\rho_0(x,y) + \sum_{\substack{(x,y)\in A \\ d(x,y_0) < d(x_0,y_0)}}\rho_0(x,y)\\
&\stackrel{\eqref{eq:proofTransportHLimpl}}{=} 2{\rho_0}(L) + \sum_{\substack{(x,y)\in X\cup Y \\ d(x_0,y) < d(x_0,y_0)}} \rho_0(x,y) + \sum_{\substack{(x,y)\in X\cup Y \\ d(x,y_0) < d(x_0,y_0)}} \rho_0(x,y) \\
&\geq 2{\rho_0}(L) + \sum_{\substack{(x,y)\in Y \\ d(x_0,y) < d(x_0,y_0)}} \rho_0(x,y) + \sum_{\substack{(x,y)\in X \\ d(x,y_0) < d(x_0,y_0)}} \rho_0(x,y) \\
&=2\rho_0(L) + \rho(X) + \rho(Y).
\end{align*}
Putting together gives
\begin{align*}
\cost(\pi) &= -2\mu_\pi(2) - \mu_\pi(1) + \mu_\pi(-1) + 2\mu_\pi(-2)\\
&\geq 2\rho_0(H) - (2\rho_0(H) + \rho_0(X) + \rho_0(Y)) + (2\rho_0(L) + \rho_0(X) + \rho_0(Y)) - 2\rho_0(L)
=0
\end{align*}
which shows that $\rho = \rho_0 + \pi$ is an optimal transport plan.

Step 3: We show $\mu_\rho(k) = 0$ whenever $|k|>1$ and $\mu_\rho(-1) \geq 2q_{\min}$.
We have
\[
\mu_\rho(2) = \mu_{\rho_0}(2)+\mu_{\pi}(2) = \rho_0(H)- \rho_0(H) = 0
\]
and similarly,
$\mu_\rho(-2) = 0$. This proves $\mu_\rho(k) = 0$ whenever $|k|>1$.

For proving $\mu_\rho(-1) \geq 2q_{\min}$, let $x' \sim x_0$ and $y' \sim y_0$ such that $d(x',y_0) = d(x_0,y') = d(x_0,y_0)-1$ and such that $d(x',y') - d(x_0,y_0) \neq -1$.
By construction, we have $\rho = 0$ on $X\cup L$ and thus, $\rho(x',y)=0$ whenever $d(x',y) - d(x_0,y_0) \neq -1$.
Similarly, $\rho(x,y')=0$ whenever $d(x,y') - d(x_0,y_0) \neq -1$.
In particular, $\rho(x',y')=0$. Thus,
\[
\mu_\rho(-1) \geq \sum_y \rho(x',y) + \sum_x \rho(x,y') = q(x_0,x') + q(y_0,y') \geq 2q_{\min}
\]
which finishes the proof of step 3.
Thus, the proof is complete.
\end{proof}

\subsection{Coupling graphs}
In this subsection, we provide a purely analytic approach to the coupling method on random walks.
The behavior of coupled random walks on $G=(V,q)$ will be encoded by coupling graphs on the set $V\times V$. The aim of this subsection is to relate coupling graphs to transport plans.

Let $G=(V,q)$ be a reversible graph with invariant measure $m$.
A graph $\widetilde G = (V\times V, \widetilde q)$ with Laplacian $\widetilde \Delta: C(V^2) \to C(V^2)$ is called a \emph{coupling graph} of $G$ if for all $f \in C(V)$,
\begin{align*}
\widetilde \Delta (f \otimes 1) = \Delta f \otimes 1 \qquad \mbox{ and } \qquad \widetilde \Delta (1 \otimes f) =   1 \otimes \Delta f. 
\end{align*}

\begin{theorem}[Coupling and transport]\label{thm:CouplingTransport}
Let $G=(V,q)$ be a reversible graph with invariant measure $m$. Let $\widetilde G= (V\times V,\widetilde q)$ be a graph on the product space $V\times V$. T.f.a.e.:
\begin{enumerate}[(i)]
\item
$\widetilde G$ is a coupling graph of $G$.
\item
For all $x_0,y_0 \in V$, the map $\rho_{x_0y_0}:V\times V \to [0,\infty)$ given by
\[
\rho_{x_0y_0}(x,y) := \widetilde q((x_0,y_0),(x,y))
\]
is a transport plan from $x_0$ to $y_0$.
\end{enumerate}
\end{theorem}
\begin{proof}
We first prove $(i)\Rightarrow (ii)$.
Let $x_0,y_0 \in V$. We aim to show that $\rho_{x_0y_0}$ is a transport plan, i.e.
\begin{align}\label{eq:proofCouplTransX}
\sum_y \rho_{x_0y_0}(x,y) = q(x_0,x)
\end{align}
for all $x \neq x_0$ and
\begin{align}\label{eq:proofCouplTransY}
\sum_x \rho_{x_0y_0}(x,y) = q(y_0,y)
\end{align}
for all $y \neq y_0$.
We first prove \eqref{eq:proofCouplTransX}. Let $x\neq x_0 \in V$.
Consider $f:= 1_x$. Then,
\begin{align*}
\sum_y \rho_{x_0y_0}(x,y)  &= \sum_y q((x_0,y_0),(x,y)) \\&= \widetilde \Delta (f\otimes 1)(x_0,y_0) \\&
 = (\Delta f \otimes 1) (x_0,y_0)
 \\&= q(x_0,x)
\end{align*}
which proves \eqref{eq:proofCouplTransX} since $x$ is arbitrary. The proof of \eqref{eq:proofCouplTransY} is analogous, and thus the proof of $(ii)$ is finished.

We finally prove $(ii) \Rightarrow (i)$. Let $x_0,y_0 \in V$.
We have
\begin{align*}
\widetilde \Delta(f \otimes 1)(x_0,y_0) &= \sum_{x,y \in V} \widetilde q((x_0,y_0),(x,y))(f(x)-f(x_0))\\
&=\sum_{x \in V \setminus \{x_0\}} (f(x)-f(x_0)) \sum_y  \widetilde q((x_0,y_0),(x,y))\\
&=\sum_{x \in V \setminus \{x_0\}} (f(x)-f(x_0))q(x_0,x)\\
&=\Delta f \otimes 1(x_0,y_0)
\end{align*}
where the third equality follows since 
$\widetilde q((x_0,y_0),(\cdot,\cdot))$ is a transport plan by $(ii)$.
Analogously, $$\widetilde \Delta(1 \otimes f)(x_0,y_0) = (1 \otimes \Delta f)(x_0,y_0).$$
This proves $(i)$ since $x_0,y_0\in V$ are chosen arbitrarily. Thus, the proof is finished.
\end{proof}

\begin{definition}
We say $\widetilde G=(V^2,\widetilde q)$ is a \emph{perfect coupling graph} of $G$ if $\widetilde G$ is a coupling graph and if for all $x_0,y_0 \in V$, the transport plan $\rho_{x_0y_0} := q((x_0,y_0),(\cdot,\cdot))$ satisfies
\begin{itemize}
\item $\rho_{x_0y_0}$ is an optimal transport plan,
\item $\mu_{\rho_{x_0y_0}}(k) = 0$ whenever $|k|>1$,
\item $\mu_{\rho_{x_0y_0}}(-1) \geq 2q_{\min}$ whenever $x_0 \neq y_0$.
\end{itemize}
\end{definition}
Putting together Theorem~\ref{thm:transportMass} and Theorem~\ref{thm:CouplingTransport} gives the following Corollary.
\begin{corollary}
Let $G=(V,q)$ a graph. Then, there exists a perfect coupling graph of $G$. 
\end{corollary}

The main idea of this article is to compare
$\widetilde P_t 1_W$ with $P_t^{\N_0}1_{\N^+}(d(\cdot,\cdot))$ where $\widetilde P_t$ is the heat semigroup on a perfect coupling graph and  $W = V^2 \setminus \diag$
and $P_t^{\N_0}$ is the heat semigroup on the non-negative integers with an absorbing state at zero. Since the vertex degree can be unbounded, we need some care to introduce the heat semigroup $\widetilde P_t$ on the potentially non-reversible coupling graph.
This is done in the next section where we also show that the semigroup of the coupling graph is compatible with tensorization.
\subsection{Heat semigroup on directed graphs and tensorization}

The aim of this subsection is to define the heat semigroup on non-reversible graphs and particularly on coupling graphs. 
The heat semigroup on reversible graphs has been investigated in e.g. \cite{keller2010unbounded,keller2012dirichlet,
keller2013volume,wojciechowski2008heat}.

We prove that the heat semigroup has all necessary properties one would expect from the reversible case. Moreover, we show that the heat semigroup of a coupling graph 
is compatible with tensorization.

Let $G=(V,q)$ be a graph.
Let $S \subset V$ be a finite subset. Let $G_S:=(V,q_S)$ with
\[
q_S(x,y) = q(x,y)1_S(x).
\]
Let $\Delta_S$ be the corresponding Laplacian.
Observe that all connected components of $G_S$ are finite and thus,
the heat equation
\[
\partial_t u_t = \Delta_S u_t
\]
with given $u_0 \in C(V)$
has the unique solution
\[
u_t = \sum_{k\geq 0} \frac{\Delta_S^k u_0}{k!}.
\]
We define
\[
P_t^S f := \sum_{k\geq 0} \frac{\Delta_S^k (f \cdot 1_S)}{k!}.
\]

\begin{theorem}[Maximum principle]\label{thm:MaxPrinciple}
Let $G=(V,q)$ be a graph.
Let $S \subset V$ be a finite subset and $T>0$.
Suppose $u_t \in C(V)$ satisfies the heat inequality
\[
\partial_t u(x) \leq \Delta u(x) \mbox{ for all } x\in S \mbox{ and all } t\in [0,T].
\]
Then,
\[
\max_{(x,t) \in S \times [0,T]} u_t(x) \leq  \sup_{(x,t) \in V \setminus S \times [0,T] \cup V \times\{0\}} u_t(x).
\]
\end{theorem}
\begin{proof}
Let $u_t^\alpha := u_t -\alpha t$ for $\alpha >0$. In order to prove the theorem, it suffices to prove that $u_t^\alpha$ as a function on $V \times [0,T]$ can not attain its maximum on $S \times (0,T]$.
Suppose $u_t^\alpha$ attains its maximum in $(x_0,t_0) \in S \times (0,T]$.
Then at $(x_0,t_0)$, 
\[
0<\alpha \leq \alpha + \partial_t u^\alpha_t 
= \partial_t u_t \leq  \Delta u_t = \Delta u_t^\alpha
\]
However $u_{t_0}^\alpha(x_0) \geq u_{t_0}^\alpha$ which implies $\Delta u_t^\alpha \leq 0$. This is a contradiction and thus, the proof is finished.
\end{proof}

We want to define $P_t f$ as a limit of $P_t^{S_n} f$ where $S_n$ are finite subsets of $V$ which eventually exhaust $V$. Using the maximum principle, we now show that this limit exists pointwise.
If $S_1 \subset S_2$ and $f \geq 0$, we have
\[
P_0^{S_2} f \geq P_0^{S_1} f
\]
and on $V\setminus S_1$, we have
\[
P_t^{S_2} f \geq 0= P_t^{S_1} f.
\]
Thus by the maximum principle,
\[
P_t^{S_2} f \geq  P_t^{S_1} f.
\]
Observe that if $f \leq 1$, then $P_S f \leq 1$ for all $S \subset V$ finite.
Let $S_1 \subset S_2 \subset \ldots$ be finite such that
\[
\bigcup_i S_i = V.
\]
Let $f \in \ell_\infty(V)$ be non-negative. Then,
\[
P_t^{S_1} f \leq P_t^{S_2} f \leq \ldots \leq \|f\|_\infty.
\]
We now define 
\[
P_t f := \sup_i P_t^{S_i} f = \lim_{i \to \infty}P_t^{S_i} f
\]
where the latter equality follows from monotonicity. Observe that the values of $P_t f$ do not depend on the exhausting sequence $(S_n)$.
For general $f \in \ell_\infty(V)$, we define 
$$P_t f := P_t f_+ - P_t f_- = \lim_{i\to \infty}P_t^{S_i} f.$$

We now show that $P_t f$ is the minimal solution to the heat equation whenever $f \geq 0$. 
\begin{theorem}[Heat semigroup as minimal solution]\label{thm:PtMinSol}
Suppose $u_t \in \ell_\infty(V)$ satisfies $u_t \geq 0$ and
\[
\partial_t u_t \geq \Delta u_t.
\]
Then,
\[
u_t \geq P_t u_0.
\]
\end{theorem}
\begin{proof}
Let $S \subset V$ be finite.
Applying the maximum principle to $P_t^S u_0 - u_t$ gives $u_t \geq P_t^S u_0$. Taking supremum over all finite $S \subset V$ proves the claim.
\end{proof}

\begin{definition}
We say a graph $G=(V,q)$ is \emph{stochastically complete} if $P_t 1 = 1$ for all $t \geq 0$.
\end{definition}

Stochastic completeness under a discrete curvature bound has been studied in \cite{hua2017stochastic,hua2017ricci,
munch2017ollivier}.

\begin{theorem}[Tensorization and stochastic completeness]\label{thm:PtTensor}
Let $G=(V,q)$ be stochastically complete with heat semigroup $P_t$. 
Let $\widetilde G =(V^2, \widetilde q)$ be a coupling graph with heat semigroup $\widetilde P_t$. Then for all $f \in \ell_\infty(V)$,
\[
\widetilde P_t(1 \otimes f) = 1 \otimes P_t f \qquad \mbox{ and } \qquad \widetilde P_t(f \otimes 1) = P_t f \otimes 1.
\]
\end{theorem}

\begin{proof}
Without loss of generality, we assume $f \geq 0$ and $\|f\|_\infty \leq 1$.
Let $S \subset V$ be finite. Let
\[
u_t^S := P_t^S f \otimes 1 - \widetilde P_t^{S\times S} (f \otimes 1).
\]
Observe that on $S\times V$,
\[
\widetilde \Delta(P_t^S f \otimes 1) = \Delta P_t^S f \otimes 1 = (\partial_t P_t^S f) \otimes 1 = \partial_t (P_t^S f \otimes 1)
\]
where the first equality holds true since $\widetilde G$ is a coupling graph. Thus, $u_t^S$ satisfies $\widetilde \Delta u_t^S = \partial_t u_t^S$ on $S^2$.
In order to prove $\widetilde P_t(f \otimes 1) = P_t f \otimes 1$ it suffices to prove that $\lim_{S \to V} u_t^S = 0$ pointwise.
Observe $u_t^S \geq 0$ by the maximum principle since $u_0 \geq 0$ and since $\widetilde P_t^{S\times S} (f \otimes 1) = 0$ on $V^2 \setminus S^2$ and since $u_t^S$ satisfies the heat equation on $S^2$.
We have $u_t \leq 1$ on $S \times V$ and $u_t = 0$ on $(V\setminus S) \times V$. We also note $u_0 = 0$ on $S^2$.
Now consider $$\phi_t^S := u_t^S + 1 \otimes P_t^S 1.$$
Observe $1 \otimes P_t^S 1 \leq 1$ on $V \times S$ and 
$1 \otimes P_t^S 1 = 0$ on $V \times (V\setminus S)$.
Thus on $(V \times (V\setminus S)) \cup ((V\setminus S) \times V) = V^2 \setminus S^2$, we have
$
\phi_t^S \leq 1.
$
Moreover, $\phi_0^S \leq 1$ on $S^2$. By maximum principle, we have $\phi_t^S \leq 1$.  
By stochastic completeness, we have
$\lim_{S\to V} 1\otimes P_t^S 1 = 1$ pointwise. Thus, $\limsup_{S \to V} u_t^S \leq 0$.
Putting together with $u_t^S \geq 0$ gives $\lim_{S\to V} u_t^S = 0$ as desired.
The equality $\widetilde P_t(1 \otimes f) = 1 \otimes P_t f$ follows analogously.
This finishes the proof.
\end{proof}

%\section{Heat semigroup on the integers}

\subsection{Comparing heat equation on coupling graphs and the integers}

In this section, we fix a graph $G=(V,q)$.
Let $G_{\N_0} = (\N_0,q_0)$ be the graph given by
\[
q_0(x,y)=\begin{cases}
2q_{\min}&: |x-y|=1 \mbox{ and } x>0,\\
0&: \mbox{else}.
\end{cases}
\]
The graph $G_{\N_0}$ can be interpreted as a birth death chain with an absorbing state at zero.
Let $P_t^{\N_0}$ be the corresponding heat semigroup.
Let $\phi_t: \N_0 \to [0,1]$ be given by
\[
\phi_t(n):=P_t^{\N_0} 1_{\N_+}(n).
\]
The aim of this subsection is to compare the heat semigroup $\widetilde P_t$ on a perfect coupling graph $\widetilde G$ to the heat semigroup $P_t^{\N_0}$ on the birth death chain given above.
Particularly, we want to show
\[
\widetilde P_t 1_W \leq \phi_t(d(\cdot,\cdot))
\]
with $W = V^2 \setminus \diag(V)$.
To this end, we first give some basic properties of the function $\phi_t$.
For geometric and analytic properties of more general birth death chains with non-negative Bakry-Emery curvature, see \cite{hua2017ricci}. 
\begin{lemma}[Heat semigoup on the integers]\label{thm:heatIntegers}
We have
\begin{enumerate}[(a)]
 \item For $r \in \N$ and $t>0$,
\[
\partial_t \phi_t(r) = 2q_{\min}(\phi_t(r+1) + \phi_t(r-1) - 2\phi_t(r)) \leq 0.
\]
\item
 $\phi_t(y) \geq \phi_t(x)$ whenever $y \geq x$.
\item For $r \in \N$ and $t>0$,
\[
\phi_t(r) \leq \frac {r} {2\sqrt {t q_{\min}}}.
\]
\end{enumerate}
\end{lemma}
\begin{proof}
We first prove $(a)$.
Let $S \subset \N_0$ be finite
We consider the graph $G_{\N_0}^S = (\N_0,q_0^S)$ given by $q_0^S(x,y) := q_0(x,y)1_S(x)$ and the corresponding Laplacian $\Delta_S^{\N_0}$ and heat semigroup $P_t^{\N_0,S}$.
By definition, we have $P_t^{\N_0,S} 1_{\N^+} \to \phi_t$ pointwise from below as $S \to V$.
Thus, it suffices to prove 
\[
\Delta_S^{\N_0} P_t^{\N_0,S} 1_{\N^+} \leq 0.
\]
Observe $\Delta_S^{\N_0} (1_S 1_{\N^+}) \leq 0$ and thus,
\[
\Delta_S^{\N_0} P_t^{\N_0,S} 1_{\N^+} =  P_t^{\N_0,S} \Delta_S^{\N_0} (1_S 1_{\N^+}) \leq 0
\]
which proves $(a)$.

We next prove $(b)$.
Suppose $\phi_t(y) < \phi_t(x)$ for some $t>0$ and some $y>x$.
By $(a)$, we have
\[
\phi_t(x + N(y-x)) \leq \phi_t(x) - N(\phi_t(x)-\phi_t(y))
\]
for $N \in \N_+$. However, the right hand side becomes negative if $N$ is large which is a contradiction since $P_t^{\N_0} g\geq 0$ for all $g\geq 0$. This proves $(b)$.

We finally prove $(c)$.
Consider $\Z$ with weights $q_{\Z}(x,y)=2q_{\min}$
if $|x-y|=1$ and $q_{\Z}(x,y)=0$ otherwise.
Let $\Delta^{\Z}$ be the corresponding Laplacian and $P_t^{\Z}$ be the corresponding heat semigroup.
Let $f(n) := \sgn(n)$ for $n \in \Z$. 
Observe
\[
P_t^{\Z}f (r) = \phi_t(r)
\]
for $r \geq 0$ since $P_t^{\Z}f (0)=0$ by symmetry which implies that  $P_t^{\Z}f$ and $\phi_t$ satisfy precisely the same heat equation
which in turn implies $P_t^{\Z}f (r) = \phi_t(r)$ by boundedness and stochastic completeness.

Since $(\Z,q_{\Z})$ has non-negative Bakry-Emery curvature (see e.g. \cite[Proposition~1.6]{lin2010ricci}), we have
\[
2t \Gamma P_t^{\Z} f \leq P_t^{\Z} f^2 - \left(P_t^{\Z}f\right)^2 \leq  1 
\]
where 
$2\Gamma g :=\Delta^{\Z}g^2 - 2g \Delta^\Z g$ and where
the first estimate follows from e.g. \cite[Theorem~3.1]{lin2015equivalent} or \cite[Theorem~1.1]{keller2018gradient}.
Thus, we have
\begin{align*}
\frac 1 {2t} \geq \Gamma P_t^{\Z} f(0) = q_{\min}\Big((P_t^{\Z} f(1)-P_t^{\Z} f(0))^2 + (P_t^{\Z} f(-1)-P_t^{\Z} f(0))^2 \Big) &= 2q_{\min} (P_t^{\Z} f(1) - P_t^{\Z} f(0))^2 \\
&= 2q_{\min} (\phi_t(1)-\phi_t(0))^2
\end{align*}
where the second equality follows from $P_t^{\Z} f(1) = -P_t^{\Z} f(-1) $ by symmetry.
Taking square root gives
\[
 |\phi_t(1) - \phi_t(0)| \leq \frac 1{2\sqrt{tq_{\min}}}.
\]
Applying assertion $(a)$ gives
\[
\phi_t(r) = \phi_t(r)-\phi_t(0) \leq r (\phi_t(1)-\phi_t(0)) \leq \frac r{2\sqrt{tq_{\min}}}
\]
for all $r \in \N$
which finishes the proof.
\end{proof}

We now show that the heat semigroup $P_t^{\N_0}$ provides a supersolution to the heat equation on a perfect coupling graph.

\begin{lemma}[Heat semigroup comparison]\label{lem:heatCompare}
Let $G=(V,q)$ be a reversible graph with non-negative Ollivier curvature.
Let $\widetilde G=(V^2,\widetilde q)$ be a perfect coupling graph with Laplacian $\widetilde \Delta$.
Let $u_t:V\times V \to [0,1]$ be given by
\[
u_t(x,y):= \phi_t(d(x,y)) = P_t^{\N_0}1_{\N_+}(d(x,y)).
\]
Then,
\[
\partial_t u_t \geq \widetilde \Delta u.
\]
\end{lemma}

\begin{proof}
First observe that
\[
\partial_t u_t(x,x)= 0 = \widetilde \Delta u_t(x,x)
\]
for all $x \in V$ since $\widetilde q((x,x),(v,w))>0$ implies $v=w$ by optimality (see Theorem~\ref{thm:transportMass}).
Now let $x \neq y\in V$ and let $r:=d(x,y)$.
Let $\rho_{xy}(v,w) = \widetilde q((x,y),(v,w))$.
We recall
\[
\mu_{\rho_{xy}}(k) = \sum_{d(v,w)-d(x,y)=k}\rho(v,w).
\]
Since $G$ is a perfect coupling graph, we have
$\mu_{\rho_{xy}}(-1)\geq 2q_{\min}$ and $\mu_{\rho_{xy}}(k)=0$ whenever $|k|>1$.
Moreover since $\rho_{xy}$ is optimal and since $\kappa(x,y) \geq 0$, we have
\[
0 \leq d(x,y)\kappa(x,y) = \cost(\rho_{xy}) = \mu_{\rho_{xy}}(-1) - \mu_{\rho_{xy}}(1).
\]
We calculate
\begin{align*}
\widetilde \Delta u_t(x,y) &= \sum_k \mu_{\rho_{xy}}(k) (\phi_t(r+k)-\phi_t(r)) \\
& = \mu_{\rho_{xy}}(-1) (\phi_t(r-1)-\phi_t(r)) + \mu_{\rho_{xy}}(1) (\phi_t(r+1)-\phi_t(r))
\\&\leq \mu_{\rho_{xy}}(-1) \Big((\phi_t(r-1)-\phi_t(r)) + (\phi_t(r+1)-\phi_t(r))\Big)\\
&= \mu_{\rho_{xy}}(-1) (\phi_t(r+1)+\phi_t(r-1)-2\phi_t(r)) \\
&\leq 2q_{\min}(\phi_t(r+1)+\phi_t(r-1)-2\phi_t(r)) \\
&=\partial_t \phi_t(r)\\
&=\partial_t u_t(x,y)
\end{align*}
where the first estimate follows since $\phi_t(r+1) \geq \phi_t(r)$ (see Lemma~\ref{lem:heatCompare}(b)) and since $\mu_{\rho_{xy}}(1) \leq \mu_{\rho_{xy}}(-1)$, and the second estimate follows since $\mu_{\rho_{xy}}(-1)\geq 2q_{\min}$ and since $\phi_t(r+1)+\phi_t(r-1)-2\phi_t(r) \leq 0$ (see Lemma~\ref{lem:heatCompare}(a)). This finishes the proof since $x,y\in V$ are chosen arbitrarily.
\end{proof}
As we have shown that $\phi_t(d(\cdot,\cdot))$ provides a super solution to the heat equation on $\widetilde G$, we can now compare it to the heat semigroup on $\widetilde G$ since the heat semigroup provides the smallest positive super solution by Theorem~\ref{thm:PtMinSol}.

\begin{theorem}[Semigroup estimate for coupling graphs]\label{thm:semigroupCouplingEstimate}
Let $G=(V,q)$ be a reversible graph with non-negative Ollivier curvature.
Let $\widetilde G=(V^2,\widetilde q)$ be a perfect coupling graph with the heat semigroup $\widetilde P_t$.
Let $W := V^2 \setminus \{(x,x):x \in V\}$.
Then,
\[
\widetilde P_t 1_W \leq \phi_t(d(\cdot,\cdot)) \leq \frac{d(\cdot,\cdot)} {2\sqrt{t {q_{\min}}}}.
\]
\end{theorem}
\begin{proof}
Let  $\widetilde \Delta$ be the Laplacian of $\widetilde G$. Let $u_t:V^2 \to [0,1]$, 
\[
u_t(x,y):= \phi_t(d(x,y)) = P_t^{\N_0}1_{\N_+}(d(x,y)).
\]

We have $\widetilde P_0 1_W = 1_W = u_0$ and $\partial_t u_t \geq \widetilde\Delta u_t$ by Lemma~\ref{lem:heatCompare}. Thus by Theorem~\ref{thm:PtMinSol},
\[
\widetilde P_t 1_W \leq u_t = \phi_t(d(\cdot,\cdot)).
\]
Moreover by Theorem~\ref{thm:heatIntegers},
$\phi_t(r) \leq \frac{r} {2\sqrt{t {q_{\min}}}}$.
Putting together with the above estimate finishes the proof.
\end{proof}

\subsection{Reverse Poincar\'e inequality}

In this subsection, we apply the semigroup estimate for coupling graphs (Theorem~\ref{thm:semigroupCouplingEstimate}) and we use the inequality $$2\|f\|_\infty 1_W \geq 1\otimes f - f\otimes 1$$ with $W=V^2 \setminus \diag(V)$ in order to prove reverse Poincar\'e inequality.

\begin{theorem}[Reverse Poinvar\'e inequality]\label{thm:NonNegOllRevPoin}
Let $G=(V,q)$ be a reversible graph with non-negative Ollivier curvature. Then for all $f \in \ell_\infty(V)$,
\[ \|\nabla P_t f\|_\infty \leq \frac {\|f\|_\infty}{\sqrt{tq_{\min}}}
.\]
\end{theorem}
\begin{proof}
Let $f\in \ell_\infty(V)$.
Let $\widetilde G=(V^2,\widetilde q)$ be a perfect coupling graph  with the corresponding heat semigroup
 $\widetilde P_t$.
Let $W := V^2 \setminus \{(x,x):x \in V\}$ and
let $x,y \in V$.
We write 
\[F:=1\otimes f - f \otimes 1.\]
Then, $F \leq 2\|f\|_\infty 1_W$.
Since $G$ has lower bounded Ollivier curvature, it is stochastically complete (see \cite[Lemma~3.6]{munch2017ollivier}) and thus by Theorem~\ref{thm:PtTensor},
\[
P_t f(y) - P_t f(x) = (1 \otimes P_t f - P_t f \otimes 1)(x,y) = \widetilde P_t F(x,y).
\]
Moreover by Theorem~\ref{thm:semigroupCouplingEstimate}, we have $\widetilde P_t 1_W \leq \frac{d(\cdot,\cdot)}{2\sqrt{tq_{\min}}}$ and thus,
\[
P_t f(y) - P_t f(x) = \widetilde P_t F(x,y) \leq 2\|f\|_\infty \widetilde P_t 1_W(x,y) \leq 2\|f\|_\infty \frac{d(x,y)}{2\sqrt{tq_{\min}}}
\]
which implies
\[
\|\nabla f\|_\infty \leq \frac{\|f\|_\infty}{\sqrt{tq_{\min}}}
\]
since $x,y \in V$ were chosen arbitrarily. This finishes the proof.
\end{proof}

\section{Applications}\label{s:Applications}
In this section, we apply the reverse Poincar\'e inequality to prove a total variation estimate,  Buser inequality, Liouville property and eigenvalue estimates.
It is well known that all these hold true in case of non-negative Bakry Emery curvature. However it was an open problem if these results also hold true in case of non-negative Ollivier curvature.
\subsection{Total variation of the heat semigroup}
In this subsection, we characterize the reverse Poincar\'e inequality by an estimate of the total variation distance of the heat semigroup of two measures. This total variation estimate can be seen as a dual version of the reverse Poincar\'e inequality.
\begin{theorem}[Total variation estimate]\label{thm:RevPoincarChar}
Let $G=(V,q)$ be a reversible graph with invariant measure $m$. Let $t>0$. T.f.a.e.:
\begin{enumerate}[(i)]
\item For all non-negative $\mu,\nu \in \ell_1(V,m)$ with $\|\mu\|_1 = \|\nu\|_1$,
\[
\|P_t \mu - P_t \nu \|_1 \leq \frac{W_1(\mu,\nu)}{\sqrt{tq_{\min}}},
\]
\item For all $f \in \ell_\infty(V)$,
\[ \|\nabla P_t f\|_\infty \leq \frac {\|f\|_\infty}{\sqrt{tq_{\min}}}
.\]
\end{enumerate}
\end{theorem}
\begin{proof}
We first prove $(i)\Rightarrow (ii)$.
Let $x,y \in V$ and let $\delta_z := 1_z/m(z)$ for $z \in V$.
For proving $(ii)$ it suffices to prove
\[
P_t f(y) - P_t f(x) \leq \frac {d(x,y)\|f\|_\infty}{\sqrt{tq_{\min}}}.
\]
We have
\begin{align*}
P_t f(y) - P_t f(x) &= \langle \delta_y - \delta_x, P_t f \rangle \\
&=\langle P_t \delta_y - P_t \delta_x, f\rangle \\
&\leq \|P_t \delta_y - P_t \delta_x\|_1 \|f\|_\infty\\
&\leq \frac{W_1(\delta_y,\delta_x)}{\sqrt{tq_{\min}}} \|f\|_\infty\\
&=\frac {d(x,y)\|f\|_\infty}{\sqrt{tq_{\min}}}
\end{align*}
where we used H\"older inequality in the first estimate and assertion $(i)$ in the second estimate. This proves $(i)\Rightarrow (ii)$.

We finally prove $(ii) \Rightarrow (i)$.
Let $f:= \sgn(P_t \mu - P_t \nu)$.
Then,
\begin{align*}
\|P_t \mu - P_t \nu \|_1 = \langle P_t \mu - P_t \nu, f \rangle = \langle \mu - \nu, P_t f \rangle \leq \frac{W_1(\mu,\nu)}{\sqrt{tq_{\min}}}
\end{align*}
where the estimate holds true
since  $\|\nabla P_t f\|_\infty \leq \frac {1}{\sqrt{tq_{\min}}}$ by assertion $(2)$, and since $$W_1(\mu,\nu) = \sup_{\|\nabla g\| \leq 1} \langle \mu - \nu, g \rangle$$ by Kantorovich duality. This proves $(ii) \Rightarrow (i)$ which finishes the proof.
\end{proof}

Putting together Theorem~\ref{thm:RevPoincarChar}
and Theorem~\ref{thm:NonNegOllRevPoin} gives the following corollary. 

\begin{corollary}\label{cor:TV}
Let $G=(V,q)$ be a reversible graph with non-negative Ollivier curvature. Let $m$ be the invariant measure.
Let $\mu,\nu \in \ell_1(V,m)$ be non-negative with $\|\mu\|_1 = \|\nu\|_1$.
Then,
\[
\|P_t \mu - P_t \nu\|_1   \leq   \frac{W_1(\mu,\nu)}{\sqrt{tq_{\min}}}.
\]
\end{corollary}

\subsection{Buser inequality}
There is a deep relationship between eigenvalues and isoperimetric constants such as Cheeger constant which we define now.
Let $G=(V,q)$ be a reversible graph. We write 
\[
|\partial A| := \sum_{x \in A} \sum_{y \in V \setminus A} w(x,y)
\]
with $w(x,y) :=m(x)q(x,y)=m(y)q(y,x)$.
We define the Cheeger constant
\[
h:= \inf_{2m(A) \leq m(V)} \frac {|\partial A|}{m(A)}.
\]
The famous Cheeger inequality states that
\[
h^2 \leq C \lambda_1. 
\]
Buser inequality is the reverse inequality, i.e.,
\[
\lambda_1 \leq C h^2
\]
 where  non-negative Ricci curvature is required.
 In case of graphs with non-negative Bakry Emery curvature, different versions of Buser inequality have been investigated in \cite{liu2014eigenvalue, liu2015curvature,liu2018buser,
 klartag2015discrete}. A version of Buser inequality under non-negative entropic curvature can be found in \cite{erbar2018poincare}.
However, if Buser inequality also holds true in case of non-negative Ollivier curvature has been an open question which is affirmatively answered in this subsection. 

\begin{lemma}
Let $G=(V,q)$ be a finite reversible graph with non-negative Ollivier curvature. Let $A \subset V$. Then,
\[
\|P_t 1_A - 1_A\|_1 \leq 2 |\partial A| \sqrt{\frac{t}{q_{\min}}}.
\]
\end{lemma}
We remark that this is precisely the same estimate as in the case of non-negative Bakry Emery curvature (see \cite[Section~I.2]{liu2014eigenvalue} where $1/q_{\min}$ is called $D_G^{nor}$).

\begin{proof}
%Let $A \subset V$.
Let $\mu := (\Delta 1_A)_+$ and $\nu := (\Delta 1_A)_-$. Then,
\[
|\partial A| = \|\mu\|_1 = \|\nu\|_1 = W_1(\mu,\nu).
\]
By Corollary~\ref{cor:TV}, we have
\[
\|\Delta P_s 1_A\|_1 = \|P_s \mu - P_s \nu\|_1 \leq \frac{W_1(\mu,\nu)}{\sqrt{s q_{\min}}} = \frac{ |\partial A|}{\sqrt{s q_{\min}}}.
\]
Integrating from $0$ to $t$ gives
\[
\|P_t 1_A - 1_A\|_1 = \left\| \int_0^t \Delta P_s 1_A ds  \right\|_1 \leq \int_0^t \|\Delta P_s 1_A\|_1 ds \leq \int_0^t \frac{ |\partial A|}{\sqrt{s q_{\min}}} ds = 2|\partial A| \sqrt{\frac {t}{q_{\min}}}.
\]
This finishes the proof.
\end{proof}

Indeed, the above lemma is the crucial step to proof Buser inequality. In particular, the proof of Buser's inequality closely follows \cite{liu2014eigenvalue}.
For convenience, we include the proof here.
\begin{theorem}[Buser inequality]\label{thm:Buser}
Let $G=(V,q)$ be a finite reversible graph with non-negative Ollivier curvature.
Let $\lambda_1$ be the smallest non-zero eigenvalue of $-\Delta$. Then,
\[
\lambda_1 \leq \frac{16 \log(2)}{q_{\min}} h^2.
\]
\end{theorem}
\begin{proof}
Let $A\subset V$ such that $2m(A) \leq m(V)$.
We write $A^c := V\setminus A$.
We have
\[
\|1_A -P_t 1_A\|_1 = 2 \langle 1_A, P_t 1_{A^c} \rangle.
\]

Let $f:=1_A/m(A) - 1_{V\setminus A}/m(V\setminus A)$.
Then, $\langle f,1 \rangle = 0$ and thus,
\[
e^{\lambda_1 t}\langle P_t f,f \rangle \leq  \|f\|_2^2 = \frac{1}{m(A)} + \frac 1 {m(A^c)} = \frac{m(V)}{m(A)m(A^c)}
\]
where the inequality follows from the spectral theorem.
Moreover, a straight forward calculation
using \[
m(A) - \langle 1_A, P_t 1_A \rangle = \langle 1_A, P_t 1_{A^c} \rangle = \langle P_t 1_A,  1_{A^c} \rangle = m(A^c) - \langle 1_{A^c}, P_t 1_{A^c} \rangle
\]
 gives
\[
\langle P_t f,f \rangle = \frac{m(V)}{m(A)m(A^c)} - \langle 1_A,  P_t 1_{A^c} \rangle \left(\frac{m(V)}{m(A)m(A^c)} \right)^2. 
\]
Hence,

\begin{align*}
e^{-\lambda_1 t}\frac{m(V)}{m(A)m(A^c)} \geq \frac{m(V)}{m(A)m(A^c)} - \frac 1 2 \|1_A -P_t 1_A\|_1   \left(\frac{m(V)}{m(A)m(A^c)} \right)^2 
\end{align*}
Dividing by $\frac{m(V)}{m(A)m(A^c)}$ and applying the above lemma gives
\[
e^{-\lambda_1 t} \geq 1 - |\partial A| \sqrt{\frac {t}{q_{\min}}}\cdot \frac{m(V)}{m(A)m(A^c)}
\]
Choosing $t_0=q_{\min} \frac{m(A)^2m(A^c)^2}{4|\partial A|^2 m(V)^2}$ gives
\[
e^{-\lambda_1 t_0} \geq \frac 1 2
\]
which implies
\[
\lambda_1 \leq \log(2)/t_0 = \frac{4\log(2)}{q_{\min}} \left(\frac{|\partial A|}{m(A)} \right)^2 \left(\frac{m(V)}{m(A^c)} \right)^2
\leq \frac{16 \log(2)}{q_{\min}} \left(\frac{|\partial A|}{m(A)} \right)^2
\]
whenever $2m(A) \leq m(V)$. Taking infimum over $A$ gives
\[
\lambda_1 \leq \frac{16 \log(2)}{q_{\min}} h^2
\]
as desired. This finishes the proof.
\end{proof}

\subsection{Liouville property}

The Liouville property states that every bounded harmonic function is constant where a function $g$ is called harmonic if $\Delta g=0$.
The Liouville property under non-negative Bakry Emery curvature and some additional assumption has been proven in \cite{hua2017liouville}.
The Liouville property under non-negative Ollivier curvature has been proven in \cite{jost2019Liouville} under the assumptions $q_{\min}>0$ and $\sup_{x\in V} \sum_{y} q(x,y)<\infty$. In this section, we show that the Liouville property still holds true when dropping the second assumption.
We also remark that the assumption $q_{\min}>0$ stands skew to the additional assumptions from \cite{hua2017liouville}.

We first give a general sufficient criterion of the Liouville property via coupling graphs.
A more general version of the below theorem in the language of random walks can be found in \cite[Theorem~2.2]{wang2012coupling}.
\begin{theorem}
Let $G=(V,q)$ be a reversible stochastically complete graph and let $\widetilde G =(V^2,\widetilde q)$ be a coupling graph of $G$.
Suppose
\[
 \widetilde P_t 1_W \to 0 \mbox{ pointwise as } t \to \infty
\]
where $W := V^2 \setminus \diag(V)$.
Then, every bounded harmonic function is constant.
\end{theorem}
In terms of random walks, the property $\widetilde P_t 1_W \to 0$ is called successful coupling which is known to imply the Liouville property, i.e., that every bounded harmonic function is constant (see e.g. \cite[Theorem~2.2]{wang2012coupling}).

\begin{proof}
Let $\phi \in\ell_\infty(V)$ be harmonic.
Due to stochastic completeness,  we have $P_t \phi = \phi$.
W.l.o.g., $\|\phi\|_\infty \leq 1$.
Let $f := 1 \otimes \phi - \phi \otimes 1$.
Then, $f \leq 1_W$ and thus,
\[
 f  = 1 \otimes P_t \phi - P_t  \phi \otimes 1 =  \widetilde P_t f \leq  \widetilde P_t 1_W  \to 0
\]
pointwise as $t \to \infty$ where the equalities follow from stochastic completeness and from Theorem~\ref{thm:PtTensor}.
Thus, $f(x,y) \leq 0$ for all $x,y \in V$ which means $\phi(y)-\phi(x) \leq 0$. This implies that $\phi$ is constant which finishes the proof.
\end{proof}
Combining the above theorem with Theorem~\ref{thm:semigroupCouplingEstimate} and the fact that non-negative Ollivier curvature implies stochastic completeness (see ~\cite[Lemma~3.6]{munch2017ollivier}) gives the following corollary.
\begin{corollary}\label{cor:Liouville}
Let $G=(V,q)$ be a reversible graph with non-negative Ollivier curvature and $q_{\min}>0$.
Then, every bounded harmonic function is constant.
\end{corollary}

\subsection{Harnack inequality}
We now give a Harnack inequality under non-negative Ollivier curvature similar to the Harnack inequality from \cite{chung2014harnack} under non-negative Bakry Emery curvature.
\begin{theorem}[Harnack inequality]\label{thm:Harnack}
Let $G=(V,q)$ be a finite reversible graph with non-negative Ollivier curvature. Let $f$ be an eigenfunction to the eigenvalue $-\lambda \neq 0$, i.e., $\Delta f = -\lambda f$. Then,
\[
\|\nabla f\|_\infty^2 \leq \frac{2e\lambda} {q_{\min}} \|f\|_\infty^2.
\]
\end{theorem}
\begin{proof}
By the reverse Poincar\'e inequality (see Theorem~\ref{thm:NonNegOllRevPoin}), we have
\[
e^{-\lambda t} \|\nabla f\|_\infty = \|\nabla P_t f\|_\infty \leq \frac{\|f\|_\infty}{\sqrt{tq_{\min}}}
\]
Choosing $t$ such that
\[
\|\nabla f\|_\infty = \frac{\sqrt{e}\|f\|_\infty}{\sqrt{tq_{\min}}}
\]
gives $1/t \leq \lambda/\log(\sqrt{e}) = 2\lambda$ and thus,
\[
\|\nabla f\|_\infty \leq {\sqrt{2e}} \sqrt{\frac{\lambda}{q_{\min}}}\|f\|_\infty.
\]
Squaring finishes the proof.
\end{proof}

\subsection{Eigenvalue estimate}

By applying the total variation estimate to the positive and negative part of an eigenfunction, we obtain an eigenvalue estimate in terms of the diameter. More precisely, $\lambda_1 \geq C/\diam(G)^2$ where $\lambda_1$ is he smallest non-zero eigenvalue of the Laplacian. An estimate of the same type under non-negative Bakry Emery curvature is given in \cite{chung2014harnack}. For a similar estimate under non-negative entropic curvature, see \cite{erbar2018poincare}.
However if this eigenvalue estimate also holds in case of non-negative Ollivier curvature has been an open question which is affirmatively answered in this subsection.

\begin{theorem}\label{thm:EigenEstimate}
Let $G=(V,q)$ be a finite reversible graph with non-negative Ollivier curvature. Let $\lambda_1$ be the smallest non-zero eigenvalue of $-\Delta$. Then,
\[
\lambda_1 \geq \log(2) \cdot \frac{q_{\min}}{\diam(G)^2}.
\]
where $\diam(G):=\max\{d(x,y):x,y \in V\}$. 
\end{theorem}
\begin{proof}
Let $\phi$ be an eigenfunction to the eigenvalue $\lambda_1 >0$. Let $\mu := \phi_+$ and $\nu:=\phi_-$. Then, $\|\mu\|_1 = \|\nu\|_1 = \frac 1 2 \|\phi\|_1>0$ and
\[
W_1(\mu,\nu) \leq \frac 1 2 \|\phi\|_1 \diam(G). 
\] 
Applying Corollary~\ref{cor:TV} gives
\[
e^{-\lambda_1 t} \|\phi\|_1=\|P_t \phi\|_1 = \|P_t \mu - P_t \nu\|_1 \leq \frac{W_1(\mu,\nu)}{\sqrt{tq_{\min}}} \leq  \frac{\diam(G)}{2\sqrt{tq_{\min}}} \|\phi\|_1
\]
for all $t>0$.
Letting $t= \diam(G)^2/{q_{\min}}$ and dividing by $\|\phi\|_1$ gives
\[
e^{-\lambda_1 \diam(G)^2/q_{\min}} \leq \frac 1 2
\]
which implies 
\[
\lambda_1 \diam(G)^2/q_{\min} \geq \log(2). 
\]
Rearranging finishes the proof.
\end{proof}

\section{Discussion}\label{s:Discussion}
We first give an informal random walk interpretation of the coupling method provided in this paper. Afterwards, we discuss the relation between the reverse Poincar\'e inequality following from non-negative Ollivier curvature and the one following from non-negative Bakry Emery curvature.
\subsection{A random walk point of view}
We now give an intuition how one can informally interpret the methods of this paper in terms of random walks.
By Theorem~\ref{thm:transportMass}, two random walks $X_t, Y_t$ on $G$ can be coupled in such a way that
\begin{itemize}
\item In each step, the distance between $X_t$ and $Y_t$ changes by at most one.
\item In each step, the probability that the distance between $X_t$ and $Y_t$ decreases by one, is at least $2q_{\min}$.
\end{itemize}
If the curvature is non-negative then, we also have that the expected distance does not increase. This means, 
\begin{itemize}
\item In each step, the probability that the distance increases is at most the probability that the distance decreases.
\end{itemize}
The distance between $X_t$ and $Y_t$ is a random variable which can informally be interpreted as a birth death chain with non-fixed transition rates, but where we do have appropriate estimates for the transition rates. These estimates allow us to compare the distance between $X_t$ and $Y_t$ to an unweighted birth death chain with fixed weights where the probability distributions of a random walk can be estimated easily. Specifically in Theorem~\ref{thm:semigroupCouplingEstimate}, we estimate 
\[
P(X_t \neq Y_t) \leq P(Z_t>0) \leq \frac{Z_0}{2\sqrt{tq_{\min}}}
\] 
where $Z_t$ is a random walk on the non-negative integers with an absorbing state $0$ which starts with $Z_0 = d(X_0,Y_0)$.

\subsection{Comparison to Bakry Emery curvature}
A reverse Poincare inequality under non-negative Bakry Emery curvature has been shown under different assumptions:
\begin{enumerate}[(A)]
\item $\Deg(x) := \sum_{y\neq x} q(x,y)$ bounded in $x \in V$, see \cite{lin2015equivalent}.
\item $\inf_{x\in V} m(x) > 0$ and there exists $\eta_k \in C_c(V)$ such that $\Gamma \eta_k \leq 1/k$ and such that $\eta_k \to 1$ pointwise from below, see \cite{gong2015properties}.
\item $q(x,y)\leq m(y)$ for all $x,y \in V$, see \cite{keller2018gradient}.
\end{enumerate}

We demonstrate that the reverse Poincar\'e inequality for non-negatve Ollivier curvature presented in this article precisely coincides with a reverse Poincar\'e inequality for non-negative Bakry Emery curvature from~\cite{lin2015equivalent,
gong2015properties,keller2018gradient,
liu2014eigenvalue,liu2015curvature}.
\begin{theorem}
Let $G=(V,q)$ be a reversible graph with non-negative Bakry Emery curvature. Suppose $G$ satisfies one of the assumptions (A), (B) or (C).
Then for non-negative $f \in C_c(V)$ and $t>0$,
\[
\|\nabla P_t f\|_\infty \leq \frac{\|f\|_\infty}{\sqrt{tq_{\min}}}.
\]
If (A) or (B) holds true, then the above inequality also holds for all $f \in \ell_\infty(V).$
\end{theorem}
\begin{proof}
The inequality
\begin{align*}
2t \Gamma P_t f \leq P_t f^2 - (P_t f)^2
\end{align*}
is implied by (A) (see \cite{lin2015equivalent}), by (B) (see \cite{gong2015properties}) and by (C) (see \cite{keller2018gradient}).
For $x\sim y$, we have
\begin{align*}
t q_{\min}|P_t f(y)-P_t f(x)|^2 \leq 2t \Gamma P_t f \leq P_t f^2 - (P_tf)^2 \leq \|f\|_\infty^2.
\end{align*}
Thus,
\[
|\nabla_{xy} P_t f| \leq \frac{\|f\|_\infty}{\sqrt{tq_{\min}}}.
\]
Taking supremum over all $x\sim y$ gives
\[
\|\nabla P_t f\|_\infty = \sup_{x\sim y} |\nabla_{xy} P_t f| \leq \frac{\|f\|_\infty}{\sqrt{tq_{\min}}}.
\]
which finishes the proof.
\end{proof}

%\subsection{Open problems}

\section*{Acknowledgments}
The author wants to thank Gabor Lippner, Mark Kempton, Christian Rose and Marzieh Eidi for fruitful discussions at Harvard University, Northeastern University and at the MPI MiS Leipzig. The author also wants to thank the MPI MiS Leipzig for financial support.

%Bobo Hua proved the Liouville property for graphs with non-negative Bakry Emery curvature \cite{hua2017liouville}. 

%\bibliographystyle{alpha}
%\bibliography{Bibliography}

\printbibliography

\end{document}